\documentclass[12pt]{amsart}
\usepackage{a4}
\usepackage{amssymb}
\usepackage{amsmath}
\usepackage{amsthm}
\usepackage{amstext}
\usepackage{amscd}
\usepackage{latexsym}
\usepackage{graphics}
\usepackage{color}

\swapnumbers
\theoremstyle{plain}
\newtheorem{thm}{Theorem}[section]
\newtheorem{q}[thm]{Question}
\newtheorem{lemma}[thm]{Lemma}
\newtheorem{lem}[thm]{Lemma}
\newtheorem{cor}[thm]{Corollary}

\newtheorem{prop}[thm]{Proposition}

\theoremstyle{definition}
\newtheorem{rem}[thm]{Remark}
\newtheorem{point}[thm]{}

\newtheorem{defn}[thm]{Definition}
\newtheorem{ex}[thm]{Example}

\newcommand{\rup}[1]{\lceil{#1}\rceil}
\newcommand{\rdown}[1]{\lfloor{#1}\rfloor}

\newcommand{\inj}{\hookrightarrow}

\newcommand{\Intersection}{\bigcap}
\newcommand{\Union}{\bigcup}
\newcommand{\intersection}{\cap}

\newcommand{\supp}{{\rm Supp}}
\newcommand{\Exceptional}{{\rm Ex}}

\newcommand{\age}{{\rm age}}

\newcommand{\Char}{{\rm char}}

\newcommand{\Gal}{{\rm Gal}}

\renewcommand{\tilde}{\widetilde}

\newcommand{\sI}{{\mathcal I}}

\newcommand{\sO}{{\mathcal O}}

\newcommand{\A}{{\mathbb A}}

\newcommand{\C}{{\mathbb C}}

\renewcommand{\P}{{\mathbb P}}
\newcommand{\Q}{{\mathbb Q}}
\newcommand{\R}{{\mathbb R}}

\newcommand{\Z}{{\mathbb Z}}

\input{xy}
\xyoption{all}

\begin{document}
\title{Degenerations of Rationally Connected Varieties}
\author{Amit Hogadi and Chenyang Xu}

\begin{abstract}
We prove that 
a degeneration of rationally connected varieties over a 
field of characteristic zero always 
contains a geometrically irreducible subvariety which is rationally connected.
\end{abstract}

\maketitle


\vspace*{6pt}

\section{Introduction}
\noindent Throughout this paper, without other mention, $k$ will denote a field of characteristic zero. A
$k$-scheme $X$ of finite type is called {\it rationally chain connected} if for an uncountable
algebraically closed field $K$ containing $k$, any two $K$-points of $X$ can be connected by a chain
of rational curves defined over $K$. A $k$-scheme $X$ is called
{\it rationally connected} if any two $K$-points of $X$ can be
connected by an irreducible rational curve defined over $K$. Thus
a rationally connected variety is always geometrically
irreducible. See \cite{kollarrc} for general background on
rational connectedness and rational chain connectedness. The
notions of rational connectedness and rational chain connectedness
are known to coincide when $X/k$ is smooth and proper. A result of
Campana and Koll\'ar-Miyaoka-Mori says that all smooth Fano
varieties are rationally connected (see\cite{Ca}, \cite{KMM} or \cite{kollarrc}). More generally, log $\Q$-Fano
varieties are rationally connected (\cite{zhang} or \cite{hm05}). 
A degeneration of a proper rationally chain
connected variety is rationally chain connected, but
a degeneration of a proper rationally connected variety may not be
rationally connected. For example a cone over an elliptic curve is
not rationally connected but it is a degeneration of smooth cubic
surfaces which are rationally connected. This paper has been
motivated by the following question raised in \cite{kollarax}.

\begin{q}{\rm (Koll\'ar)}\label{mainq'}
Suppose we have a family of smooth proper rationally connected varieties
over a field $k$. Does the degeneration $X_0$ contain an rationally connected (in particular geometrically irreducible) subvariety ?
\end{q}

\noindent In \cite{kollarax} it was proved that 
a degeneration of Fano varieties over a smooth curve always contains a geometrically irreducible subvariety.
A recent paper by Jason Starr \cite{starr}
shows that if the field $k$ (of any characteristic) contains an algebraically closed field, then
a degeneration of smooth proper seperably rationally connected varieties contains
a geometrically irreducible subvariety. In fact, A.J.de Jong proved, under the same assumption as in \cite{starr}, if $\Char(k)=0$, then the degeneration contains a rationally connected variety. As observed in \cite{kollarax}, the
existence of a geometrically irreducible subvariety can fail even for a degeneration of elliptic curves. \\

\noindent  The following result, which is the main result of this paper, gives an affirmative answer to (\ref{mainq'}). 

\begin{thm}\label{deg-ran}
Let $\pi:X \to C$ be a dominant proper morphism of $k$ varieties such that
\begin{enumerate}
\item $C$ is $klt$, and
\item The generic fibre of $\pi$ is a rationally connected variety.
\end{enumerate}
\noindent Then for any point $0 \in C$, $X_0=\pi^{-1}(0)$ contains a subvariety (defined over $k(0)$) which is geometrically irreducible and rationally connected.
\end{thm}

\noindent Our approach to (\ref{mainq'}) uses the Minimal Model Program, particularly as developed in \cite{mmpgen}. After running the minimal model program for a suitable birational model for X we obtain a Fano fibration. This is then used to reduce (\ref{deg-ran}) to the case when the generic fibre of $\pi$ is 
a log Fano variety (\ref{mainthm}). We first prove (\ref{mainthm}) in the special case when $C$ is a smooth curve and $0\in C$ is a $k$-point. The general case of arbitrary $klt$ base $C$ is then reduced to this case using the following result. 

\begin{thm}\label{exc.thm}
Let $(X,\Delta)$ be a $klt$ pair with 
$\Delta$ effective. Let $\pi:X'\to X$ be a proper birational morphism. Then
for any point $p\in X$, there exists a rationally connected subvariety 
of $\pi^{-1}(p)$ (defined over $k(p)$).
\end{thm}

\noindent The following definition was suggested by J. Koll\'ar.
\begin{defn}
A field $k$ is called {\it geometrically-$C_1$} if every smooth proper rationally connected
variety over $k$ has a $k$-point.
\end{defn}

\noindent Following are some of the consequences of (\ref{deg-ran}). 
\begin{cor}\label{fanocor}
Every geometrically-$C_1$ field (of characteristic 0) is $C_1$.
\end{cor}
\begin{proof}
To show  that $k$ is $C_1$  we have to show that every hypersurface $H_d$ of degree $d$ in $\P^n_k$
has a $k$-point if $n \geq d$. If $H_d$ is smooth, it is Fano and hence
rationally connected and has a $k$-point by the given hypothesis. Even if $H_d$ is not smooth it
can be expressed as a degeneration of smooth Fano hypersurfaces and has a a rationally
connected subvariety $Z$ by (\ref{deg-ran}). Then $Z$ has a $k$-point since
its resolution $Z'$ is smooth proper and rationally connected and has a $k$-point
by hypothesis.
\end{proof}

\begin{cor}\label{surfacecor}
Let $k$ be a $C_1$ field. Then every degeneration
of a proper rationally connected surface defined over $k$ has a $k$-point.
\end{cor}

\begin{rem} (In this remark, $k$ is not necessary of characteristic 0).
The converse of (\ref{fanocor}) is a well known problem.
It is known that a smooth proper rationally connected variety $X/k$ has a $k$-point if
\begin{enumerate}
\item $k$ is any $C_1$ field and $\dim(X)\leq 2$ (see for example (\cite{kollarrc}, {\rm IV.6.8})),
\item $k$ is a field of transcendence degree one over an algebraically closed field
of characteristic zero \cite{ghs} (in the positive characteristic case one needs to
assume $X$ is separably rationally connected \cite{ds03}), or
\item $k$ is a finite field \cite{He}.
\end{enumerate}

\noindent However, the converse of (\ref{fanocor}) is unknown even for specific $C_1$ fields like the 
maximal unramified extension of local fields.

\noindent In view of (\ref{deg-ran}), if a field $k$ is geometrically-$C_1$,
then not only do smooth proper rationally connected varieties over $k$ have 
$k$-points, but one also gets $k$-points on the degenerations of such varieties.
In fact, besides (\ref{surfacecor}), the following was already known.
\begin{enumerate}
\item A degeneration of
smooth rationally connected varieties over a finite field has
a rational point (\cite{nr}, \cite{ex07}).
\item A degeneration of seperably rationally connected variety over a field of
transcendence degree one over an algebraically closed field has a rational point
(this result has been mentioned in \cite{starr} and can be derived by using 
the properness of Kontsevich's moduli space of stable maps).
\end{enumerate}
\end{rem}

\noindent A stronger version of (\ref{mainq'}) was asked in \cite{kollarax}.

\begin{q}{\rm (Koll\'ar)}\label{mainq}
For a flat family of proper varieties over a smooth $k$-variety whose general fiber is
smooth and rationally connected, does every simple normal crossing ($snc$)
fiber $X_0$ have an irreducible component (defined over $k$) which is rationally connected (and
thus geometrically irreducible)?
\end{q}

\begin{rem}
Unlike (\ref{mainq'}), the answer to this question is unknown even if $k$ is algebraically closed.
Over an algebraically closed field ${\overline{k}}$ one can ask the following
more general question: Does a $snc$ rationally chain connected variety $X_0/{\overline{k}}$
always contain an irreducible component which is rationally connected? The answer to
this question is indeed true if $\dim(X_0)\leq 2$ but is false in higher dimensions (see (\ref{counterex})).
\end{rem}

\noindent We show that (\ref{mainq}) has an affirmative answer in the case that the generic fiber satisfies $\dim(X_t) \le 3$.   
\begin{thm}\label{surface}
Let $\pi:X/k \to C/k$ be a dominant proper morphism
such that
\begin{enumerate}
\item $C$ is a smooth curve with a $k$-point $0$,
\item $\pi$ is smooth outside $0$,
\item The generic fibre of $\pi$ is geometrically rationally connected and $\dim(X_t)\le 3$, and
\item $X$ is smooth, and $X_0= \pi^{-1}(0)$ is $snc$ (over $k$).
\end{enumerate}
\noindent Then there exists an irreducible component of $X_0/k$ which is rationally connected.
\end{thm}

\noindent In Section \ref{section:fano} we prove Theorem (\ref{deg-ran}) in the case of degeneration of log Fano varieties, assuming Theorem (\ref{exc.thm}). In Section \ref{section:mmp} we prove Theorem (\ref{deg-ran}) using the results of the minimal model program as proved in (\cite{mmpgen}). The proof uses results of Section \ref{section:fano} and some results on Fano fibrations which are proved in Section \ref{section:surface}. The stronger result, Theorem (\ref{surface}), for degeneration of rationally connected threefolds (or lower dimension) is proved in Section \ref{section:surface}. Finally, in Section \ref{section:exc} we prove Theorem (\ref{exc.thm}).  \\

\noindent{\bf Acknowledgement}: We are indebted to our advisor, J\'anos Koll\'ar, for many
useful discussions and crucial suggestions. The relevance of \cite{hm05} and
minimal model program was pointed out to us by him. We would also like to thank
Tommaso de Fernex and H\'el$\grave{\rm e}$ne Esnault for useful suggestions and Yuri Prokhorov for informing us his example (\ref{counter}).
The second author would like to thank Jason Starr for 
useful discussions. We are also grateful to the anonymous referee for enormous suggestions and corrections.

\section{Summary of known results}

\noindent In this section we first recall the related results from \cite{hm05} and \cite{kollarax}. In our proof, their ideas are heavily used. Then we state another main ingredient which is from the minimal model program theory (cf. \cite{kollarmori}) and proved in \cite{mmpgen}. \\ \\ \noindent The following is the main theorem of \cite{hm05}, stated in a slightly less generality.

\begin{thm}[\cite{hm05}]
  Let $(X,\Delta)$ be a $klt$ pair. Let $g:Y\to X$ be a birational morphism such that $-K_X$ is relatively big and $\sO_X(-m(K_X+\Delta))$ is relatively generated for some $m>0$. Then every fiber of the the composite morphism $\pi:Y\to S$ is rationally chain connected. 
\end{thm}

\noindent Specialising to the case when $g$ is the identity morphism, $(X,\Delta)$ is rationally chain connected and $Y\to X$ is a resolution of singularities, one gets the following important corollary.

\begin{cor}[\cite{hm05}]
Let $(X,\Delta)$ be a $klt$ pair which is rationally chain connected. Then $X$ is rationally connected.
\end{cor}
\noindent Recall that there exist singular varieties (for example the cone over an elliptic curve) which are rationally chain connected but not rationally connected. However, the above corollary says that such examples necessarily have singularities worse than $klt$. \\ 

\noindent In Section \ref{section:fano}, we consider the degeneration of Fano varieties (see (\ref{mainthm})), which has also been considered in \cite{hm05}. Part of (\cite{hm05}, 5.1) guarantees the existence of an irreducible component which is rationally chain connected modulo the intersection of some other components.
But this component may be neither geometrically irreducible nor rationally connected.
Although the statement of (\cite{hm05}, 5.1) does not imply (\ref{mainthm}),
a modification of methods used in (\cite{hm05}) will be enough to prove it. \\

\noindent In their paper, They prove the following technical result (\cite{hm05}, 4.1), which we only state in a weaker form we need.

\begin{prop}[\cite{hm05}, 4.1]\label{criterion}
Let $(F,\Delta)$ be a projective log pair, and let $t : F\dashrightarrow Z$ be a dominant rational map, respectively, where $F$ and $Z$ are projective,
with the following properties:
\begin{enumerate}
\item the locus of log canonical singularities of $(K_F+\Delta)$ does not dominate $Z$;
\item $F+\Delta$ has Kodaira dimension at least zero on the general fibre of $F\dashrightarrow Z$ (i.e.,
if $g : F'\to F$ resolves the indeterminacy of $F\dashrightarrow Z$, then $g^*(K_F + \Delta)$ has
Kodaira dimension at least zero on the general fibre of the induced morphism $F' \to Z$;
\item $K_F+\Delta$ has Kodaira dimension at most zero; and
\item there is an ample divisor $A$ on $F$ such that $A\le\Delta$.
\end{enumerate}
Then either $Z$ is a point, or it is uniruled.
\end{prop}
 
\noindent The main application of this proposition is when $t$ is the rational map from $F$ to its maximal rationally connected fiberation (cf. \cite{kollarrc}, IV Theorem 5.4). In fact, with the help of the following lemma, we obtain a quite useful criterion to show the rational connectedness of varieties.

\begin{lemma}[\cite{hm05}, 4.2(i)]\label{MRC}
Let $F$ be a normal variety. $F$ is rationally connected if and only if, for every nonconstant dominant
rational map $t: F \dashrightarrow Z$, $Z$ is uniruled.
\end{lemma}
\noindent Next we recall the main result of \cite{kollarax}.

\begin{thm}[\cite{kollarax}] Let $k$ be a field of characteristic zero and $C$ be a smooth curve, $Z$ an irreducible projective variety and $g:Z\to C$ a morphism. Assume that the general fibers $F_{gen}$ are
\begin{enumerate}
\item smooth,
\item geometrically connected, and
\item Fano.
\end{enumerate}
Then for every point $c$ of $C$, the fiber $g^{-1}(c)$ contains a $k(c)$-subvariety which is geometrically irreducible. 
\end{thm}

\noindent In our proof, as a special case, and also as a middle step, we improve Koll\'ar's result by showing that the subvariety which he looks at is not only geometrically irreducible but also rationally connected (see (\ref{specialcase})). Then Theorem(\ref{deg-ran}) generalizes the above result to the case where the general fibers are rationally connected and $C$ is any $klt$ base. Moreover (\ref{deg-ran}) guarantees the existence of a subvariety which is in fact rationally connected.\\  

\noindent Finally, in the very important paper \cite{mmpgen}. A large part of the minimal model program is proved. In particular, we state the following theorem which is a special case of their results. This will be the main tool in the proof of (\ref{deg-ran}).

\begin{thm}[\cite{mmpgen}]\label{mmp}
Let $(X,\Delta)$ be a $klt$ pair, projective over a base $S$. Assume $K_X+\Delta$ is not pseudoeffective. Then one can run a relative minimal model program for $(X,\Delta)$ to obtain a Fano fibration, i.e. there exists a sequence of maps
$$ X=X_0\stackrel{\phi_1}{\dashrightarrow}X_1 \stackrel{\phi_2}{\dashrightarrow} \cdots \stackrel{\phi_{r}}{\dashrightarrow}X_r  \stackrel{h}{\to} Y $$
where each $\phi_i$ is either a divisorial contraction or a flip and $h$ is a Fano fibration of relative Picard number one.
\end{thm}

\noindent The important case in which the above hypothesis is satisfied is when $X$ is a smooth rationally connected variety and $\Delta$ is empty. Indeed, in this case, the canonical divisor $K_X$ has negative intersection with any free rational curve and hence it cannot be pseudoeffective.

\section{Degenerations of Fano Varieties}\label{section:fano}

\noindent In this section we prove (\ref{deg-ran}) in the case when the generic fibre of $\pi:X\to C$ is a log Fano variety. Recall that a proper variety $X/k$ is called {\it log Fano} (or {\it log $\Q$-Fano}) if there exists an effective divisor $D$ on $X$ such that $(X,D)$ is $klt$ and $-(K_X+D)$ is nef and big. It is known that a log Fano variety is always rationally connected (cf. \cite{zhang}, \cite{hm05} ). 

\begin{thm}\label{mainthm} Let $\pi:X \to C$ be a dominant proper morphism of $k$ varieties such that
\begin{enumerate}
\item $C$ is $klt$, and
\item The generic fibre of $\pi$ is a log Fano variety.
\end{enumerate}
\noindent Then for any point $0 \in C$, $X_0=\pi^{-1}(0)$ contains a subvariety (defined over $k(0)$) which is geometrically irreducible and rationally connected (over $k(0)$). In another words, Theorem (\ref{deg-ran}) is true in the case when the generic fibre of $\pi$ is a log Fano variety.
\end{thm}

\noindent The extension theorem (\cite{hm06}, 3.17) lies at the heart of the proof. We have to first introduce the terminology there.

\begin{defn}{\rm (\cite{hm06}, $3.7$)} Let $\pi:(X,\Delta)\to S$ be a relative log pair. 
A Cartier divisor $D$ on $X$ will be called {\it $\pi$-transverse} to $(X,\Delta)$
if the natural map 
$$\pi^*\pi_*(\sO_X(D))\to \sO_X(D)$$ is sujective at the generic point of every log canonical centre
of $(X,\Delta)$. Obviously, the notion only depends on the divisor class of $D$. Similarly, a $\Q$-divisor $D$ will be called {\it $\pi$-$\Q$-transverse} to
$(X,\Delta)$ if $mD$ is $\pi$-transverse
for some $m>0$. 
\end{defn}

\noindent The proof of (\ref{mainthm}) relies on the following extension theorem
proved in \cite{hm06}. The statement given below is less general
than the one in \cite{hm06}.

\begin{thm}{\rm (\cite{hm06}, 3.17)}\label{extensionthm}
Let $\pi:X\to S$ be a projective morphism where $X/k$ is a smooth
variety of dimension $n$.
Let $F\subset X$ be a smooth divisor.
Let $B$ be effective $\Q$-divisor on $X$ such that
\begin{enumerate}
\item $B+F$ is $snc$, 
\item the support of $B$ and $F$ do not have common component, 
\item $\rdown{B}=0$,
\item $K_F+ B_{|F}$ is effective, and
\item $K_X +B+F$ is $\pi$-$\Q$-transverse to $(X,\rup{B}+F)$.
\end{enumerate}
For every sufficiently $\pi$-ample divisor $H$ on $X$ and for every sufficiently
divisible positive integer $m$, the image of the morphism
$$ \pi_*(\sO_X( m(K_X+B+F)+(n+2)H))\to \pi_*(\sO_F(m(K_F+B_{|F})+(n+2)H_{|F}))$$
contains the image of
$$ \pi_*(\sO_F(m(K_F+B_{|F})+H_{|F}))\stackrel{(n+1)H_{|F}}{\longrightarrow} \pi_*(\sO_F(m(K_F+B_{|F})+(n+2)H_{|F})).$$
\end{thm}

\noindent In the situation of our interest, condition (4) above
will be checked as follows. There will be an effective $\Q$-divisor $A$ such that
\begin{enumerate}
\item[(i)]  $K_X+F+B\sim A$, and
\item[(ii)] $A$ and $B+F$ have no common components and $A+B+F$ is $snc$.
\end{enumerate}
Since $\rup{B}+F$ is a $snc$ divisor where each irreducible
component occurs with coefficient one, the log canonical centres of $(X,\rup{B}+F)$
are exactly all possible intersections of the components of $\rup{B}+F$. But since
$A+B+F$ is $snc$ and $A$ and $B+F$ have no common components,
$\sO_X(A)$ is generated by its global sections at the generic point of every log
canonical center of $(X,\rup{B}+F)$.

\noindent As we mentioned before, the proof is given first when $C$ is a smooth curve. 

\begin{prop} \label{specialcase}
When $C$ is a smooth curve, Theorem (\ref{mainthm}) is true.
\end{prop}

\begin{proof} By extending the base field, we can assume $0\in C$ is a $k$-point. By assumption, the generic fiber $X_g$ of $\pi:X\to C$ is a log Fano variety. Thus there exists an effective divisor $D_g$ on $X_g$ such that $(X_g,D_g)$ is $klt$ and $-(K_{X_g}+D_g)$ is nef and big. Thus there exists an ample divisor $H_g$ and an effective divisor $E_g$ on $X_g$ such that $-(K_{X_g}+D_g) \sim H_g+E_g$. Furthermore, since $-(K_{X_g}+D_g)$ is nef and big, $E_g$ and $H_g$ can be chosen such that $(X_g,D_g+E_g)$ is $klt$ (cf. \cite{kollarmori} (2.61)). By replacing $D_g$ by $D_g+E_g$ we may assume $-(K_{X_g}+D_{g})\sim H_g$ is ample. Let $D$ be an effective divisor on $X$ such that $D_{|X_g}=D_g$. By shrinking $C$ we may assume that $C$ is affine and $D$ does not contain any vertical divisor outside $0$.  

\noindent After successively blowing up $X$ with suitable centers and shrinking $C$ one can find a birational morphism $h: Y \to X$ such that

\begin{enumerate}
\item $Y$ is smooth, and the natural projection $\pi'(=\pi\circ h):Y\to C$ is smooth outside $0$, and
\item $\supp (h^{-1}D + Y_0) +\Exceptional (h)$ is $snc$ (over $k$) in $Y$. Here $Y_0=h^{-1}(0)$.
\end{enumerate} 
\noindent Then there exist effective $\Q$-Cartier divisors $G$ and $R$ on $Y$ having no common components
in their support and such that $K_{Y/X} +  G- R= h^*D$. Moreover, $G$ and $R$ satisfy
\begin{enumerate}
\item Any irreducible component of $G$ which is not supported in $Y_0$ appears in $G$ with coefficient strictly less than one (this is a consequence of the fact that $(X_g,D_g)$ is $klt$), and 
\item $R$ is $h$-exceptional.
\end{enumerate}

\noindent  Note that the restriction of  $-(K_Y+G-R)$ to the generic fiber $Y_g$ of $Y$ over $C$ is the pull back of an ample divisor on $X_g$. Thus for some sufficiently small $h$-exceptional effective divisor $E$ on $Y$, $-(K_Y+G-R+E)$ is ample on the generic fibre. If necessary blow up $Y$ along smooth centers in the fiber $Y_0$ so that there exists a vertical $\Q$-Cartier divisor $\tilde\Delta$, supported in $Y_0$, such that $-(K_Y+G-R+E+\tilde\Delta)$ is relatively ample over $C$ (see \cite{kollarax}, page $5$).

\noindent Let $\tilde{\Delta} = \sum d_i \tilde{\Delta}_i$
where the divisors $\tilde{\Delta}_i$ are the irreducible components of the support $\tilde{\Delta}$. The only irreducible components of $G-R+E+\tilde{\Delta}$ which appear with coefficient $\geq 1$ are supported in $Y_0$. Thus by perturbing the coefficients $d_i$ slightly (by small rational numbers) and by changing $\tilde{\Delta}\to \tilde{\Delta}+rY_0$ for a suitable $r\in \Q$ we may assume 
\begin{enumerate}
\item $-(K_Y+G-R+E+\tilde{\Delta})$ is $\pi$-ample, $\supp(\tilde{\Delta})\subset X_0$,
\item There exists an irreducible divisor $F$ (over $k$) in $Y$, an effective $\Q$-Cartier divisors $A$, and a boundary $\tilde{B}$ in $Y$ such that $G-R+E+\tilde{\Delta}=F+\tilde{B}-A$ and such that $F$ and the supports of $A$ and $\tilde{B}$ have no pairwise-common components, and
\item There exist effective $\Q$-Cartier divisors $A_0$ and $R_0$ in $Y$ such that
$A=A_0+R_0$, $\supp(A_0) \subset Y_0$ and $\supp(R_0)$ is in the
exceptional set of $h$ and contains no component of $Y_0$. 
\end{enumerate}
The necessary perturbations as well as the proof of the existence of $F$,  $\tilde{B}$ and $A$ can be found
in \cite{kollarax}. Notice that $A+\tilde{B}+F$ is $snc$ since it is contained in $\supp (h^{-1}D + Y_0 +\Exceptional (h))$. In (\cite{kollarax}, Remark 8) it was proved that $F$ is geometrically irreducible. If $A$ were equal to $0$ then (\cite{hm05}, 5.1) would imply that $F$ is rationally connected. But in general $A$ can be nonzero. Fortunately, as we will show below, we are saved in this situation by the fact (3).\\

\noindent  Let $m$ be a sufficiently large and divisible integer and let $L_m$ be the Cartier divisor of a general section of $\sO_X(m(-K_Y-F-\tilde{B}+A))$. Let $L= \frac{1}{m}L_m$. Let $B= \tilde{B} + L$ and $\Delta = F+B-A$. Thus we are finally in the following situation 
\begin{enumerate}
\item[(i)] $-(K_Y+\Delta)\sim_{\pi'} 0$,
\item[(ii)] $\Delta = F+B-A$ where  $A+B+F$ is $snc$, 
\item [(iii)] $F$ is geometrically irreducible $\supp( F) \subset \supp (Y_0)$,
\item [(iv)] $A$, $B$ are effective $\Q$-divisors and $\rdown{B}=0$,
\item[(v)] $A=A_0+R_0$ such that $\supp (A_0)\subset \supp(Y_0)$, $R_0$ is $h$-exceptional, and
\item[(vi)] There exists a $\pi'$-ample divisor $L$ such that $L \leq B$.
\end{enumerate}

\noindent To apply (\ref{MRC}), we now make the following two claims.
\begin{enumerate}
\item[]\hspace{-6mm}Claim(1): Let $p:F'\to F$ be a proper birational morphism
such that the induced rational map $f':F'\to Z$ is a morphism. Then $p^*(K_{F}+B_{|F})$
has Kodaira dimension at least zero on the general fibre of $f'$.
\item[]\hspace{-6mm}Claim(2): $\kappa(K_F+B_{|F})=0$.
\end{enumerate}
Notice that $(F,B_{|F})$ is $klt$ since $B+F$ is $snc$, $B$ and $F$ do not have common components
and $\rdown{B}=0$. Also the support of $B$ contains the support of an
ample divisor $L$. Hence once we prove the above two claims it will
follow that the hypotheses of (\ref{criterion}) are implied by our hypotheses. In fact, the first hypothesis of (\ref{criterion}) follows from the conclusion that $(F,B|_F )$ is klt, Claim 1 implies the second hypothesis,
Claim 2 implies the third hypothesis and the fourth hypothesis follows easily from the fact that $B+F$ is $snc$. So we can conclude that $Z$ is uniruled, thus proving the rational connectedness of $F$. \\

\noindent By restricting $K_Y + F + B \sim_{\pi '} A$ on $F$, we have $K_F+ B_{|F} \sim A_F $. But
$A$ and $F$ have no common components and hence $A_F$ is effective. Thus $p^*(K_F + B_{|F})$
is effective. Hence it is also effective when restricted to the general fibre of $f'$. This
proves Claim(1).
\noindent To prove Claim(2) we use (\ref{extensionthm}). Let $M$
be a sufficiently $\pi '$-ample divisor. To prove
$\kappa(K_F+B_{|F})=0$, it is enough to show that for $m$ large
and divisible, $h^0(F,m(K_F+B_{|F})+M)= length \left(\pi
'_*(\sO_F(m(K_F+B_{|F})+M_{|F})) \right)$ is bounded independent
of $m$. By (\ref{extensionthm}) it is enough to bound the rank of
the torsion free sheaf $\pi'_*(\sO_Y( m(K_Y+B+F)+(n+2)M))$.
$K_X+B+F\sim_{\pi '}A$ and hence we have to bound the rank of $\pi
'_*(\sO_Y(mA+(n+2)M))$. Since $A=A_0+R_0$ where $R_0$ is effective and $h$-exceptional, 
this sheaf is isomorphic to a subsheaf of $\pi_*(\sO_X(mh_*(A_0) + h_*((n+2)M)))$. But since 
$\supp(h_*(A_0))\subset \supp(X_0)$, $\pi_*(\sO_X(mh_*(A_0) + h_*((n+2)M)))$ is isomorphic
to $\pi_*(\sO_X( h_*((n+2)M)))$ on $C-\{0\}$ and hence it has
constant rank.
\end{proof}

\begin{proof} [ Proof of (\ref{mainthm}) in General Case] To prove Theorem (\ref{mainthm}) for any $klt$ variety $C$ and an arbitrary point $0\in C$, we need Theorem (\ref{exc.thm}) in the following way.
Let $h:C'\to C$ be a resolution of $C$, then by (\ref{exc.thm}), there exists a rationally connected subvariety $Z \subset h^{-1}(0)$. Let $L$ be the function field of $Z$. Since $Z$ is itself rationally connected, to find a rationally connected subvariety of $X_0$, it suffices to find a rationally connected subvariety of $X_0\times_kL$ by \cite{ghs}.
Thus by base extending to $L$ we may assume there is a $k$-point $0' \in C'$ such that $h(0')=0$. Now let $D\subset C'$ be a general smooth curve passing through $0'$. Note that it is possible to find such a curve since $C'$ is smooth. Let $X_D= X\times_CD$. Since $D$ is general, the generic fibre of $X_D\to D$ is a log Fano variety. Thus (\ref{mainthm}) now follows from the above case.
\end{proof}


\section{Degenerations of Rationally Connected Varieties and the 
Minimal Model Program}\label{section:mmp}

\noindent In this section we prove Theorem(\ref{deg-ran}). It suffices to prove Theorem(\ref{deg-ran}) after replacing $X$ by any other higher birational model $X'\to X$. Thus by resolution of singularities, we may assume $X$ is smooth, and $X_0=\pi^{-1}(0)$ is $snc$ (over $k$). 
The proof of Theorem(\ref{deg-ran}) will proceed by studying the relative minimal 
model program for the pair $(X,\Delta)$ where $\Delta$ denotes the sum of all 
irreducible components of $X_0$, each appearing with coefficient one. In this situation
$X$ is smooth and $\Delta$ is $snc$. 

\begin{lem}\label{dltklt}
Let $\phi:(X,\Delta)\dashrightarrow (X',\Delta')$
be a birational morphism obtained by running the relative minimal model program
for $(X,\Delta)$ (i.e. $\phi$ is a composition of extremal divisorial contractions and
flips). Then
\begin{enumerate}
\item[(i)] Every irreducible component $\Delta'_j$ of $\Delta'(=\sum_{j=1}^m\Delta'_j)$ is of the
form $\phi_*\Delta_i$ for some $i$, and
\item[(ii)] Every $k$-irreducible component of $\Delta'$ is $klt$.
In particular, after running the relative minimal model program, 
the $\Gal(\bar{k}/k)$-conjugate components of $\Delta\times_k\overline{k}$ do not intersect.
\end{enumerate}
\end{lem}

\begin{ex}
Let $Y=(x^2+y^2+tz^2=0)\subset \P^2_{\R} \times \A^1_{\R}$ where $t$ is the
coordinate on the affine line $C:=\A^1_{\R}$. Let $X$ be the blow up of
$Y$ at $([0,0,1],0)$, and denote the exceptional divisor to be $E$.
Then the projection $f: X \to C$ gives us a family which is
defined over $\R$, whose generic fibres are smooth conics. Look at
the point $t=0$ on $C$, the fibre $X_0$ consists of $E$ and the
birational transform of $Y_0=(x^2+y^2=0)$, which has 2
components $E_1, E_2$ over $\C$, and are conjugate to each other
under the Galois action of $\Gal(\C/\R)$. So $NE(X/C)(\R)$ is
generated by $E$ and $E_1+E_2$. We have $(K_X+E+E_1+E_2)\cdot E=0$
and $(K_X+E+E_1+E_2)\cdot (E_1+E_2)=-2$, so to run the relative
program for $f:(X,E+E_1+E_2) \to C$ over $\R$ , we can only
contract $E_1+E_2$.
\end{ex}

\begin{proof}[Proof of (\ref{dltklt})]
$(i)$ For arbitrary  $\phi$
as in Lemma 3.1, it is a morphism away from codimension 2 subsets of both the
domain and target, hence (i) is obvious. \\
\noindent $(ii)$ We devide the argument into several steps.\\
\noindent$Step(1)$ Let $\Delta'= \sum_{i=1}^r\Delta_i'$ where $\Delta'_i= \phi_*\Delta_i$ and
$\phi_*\Delta_i = 0 \ \forall \   r+1 \leq i \leq m$.
By Inversion of Adjunction, to prove $\Delta_1'$ is $klt$ it is
enough to prove $(X',\Delta'_1)$ is $plt$ (\cite{kollarmori}, $5.4$). For $I\subset \{1,...,m\}$ let $\Delta_I= \Intersection_{i\in I}\Delta_{i}$.
Since $X$ is smooth and $\Delta \subset X$ is $snc$,
$(X,\Delta)$ is log canonical and the only log canonical centers of $(X,\Delta)$ are
connected components of closed subsets $\Delta_I$.
Because during each step of the minimal model program discrepancy never decreases 
(\cite{kollarmori}, $3.38$), $(X',\Delta')$ is also log canonical. \\
\noindent $Step(2)$ Let $E$ be an exceptional divisor of $X'$
with $a(E,X',\Delta')=-1$. We need to show $a(E,X',\Delta_1')>-1$.
We claim that $center_{X}E = \Delta_I$ for some $I$. Applying the nondecreasing of the discrepancy again,  we conclude $a(E,X,\Delta)=-1$.
But as remarked earlier the only log canonical centers of $(X,\Delta)$ are $\Delta_I$.
Hence $center_XE=\Delta_I$ for some $I$. \\
\noindent $Step(3)$ We claim that for every $i\in I$, $\phi_*\Delta_i\neq 0$, i.e.,
none of the components of $\Delta$ which contain the generic point of $center_XE$ can be
contracted to a lower dimension subvariety by $\phi$. This is a consequence of
(\cite{kollarmori}, $3.38$) since if $\phi$ is not an isomorphism at the generic point
of $center_XE$ then $$a(E,X',\Delta')>a(E,X,\Delta)\geq -1$$ which is contradictory to
our assumption. Since we have assumed $E$ is an exceptional divisor of $X'$,
this also shows that $center_XE$ is not a divisor. Thus $I$ contains at least two elements.
Without loss of generality we assume $2 \in I$.\\
\noindent $Step(4)$ Now we have $a(E,X',\Delta_1'+\Delta_2')\geq a(E,X',\Delta') = -1$.
But both $\Delta_1'$ and $\Delta_2'$ contain $center_{X'}E$ hence
$$ a(E,X',\Delta_1') > a(E,X',\Delta_1'+\Delta_2') \geq -1$$
\end{proof}

\noindent The following lemma is useful in proving Theorem(\ref{deg-ran}) because it allows us to replace $X$ by any other $klt$ birational model.

\begin{lemma}\label{birational}
Let $f:X'\dasharrow X$ be a birational map of varieties over $C$ such that both $X'$ and $X$ are $klt$. Then for $0\in C$, $X'_0$ has a rationally connected subvariety if and only if the same is true for $X_0$. 
\end{lemma}
\begin{proof}Suppose $X_0$ contains a rationally connected subvariety $Z$. It is enough to prove $X_0'$ contains a rationally connected subvariety after replacing $X'$ by a higher birational model. Thus we may assume $f$ is in fact a morphism. Let $k(Z)$ be the generic point of $Z$. By Theorem(\ref{exc.thm}), there exists a rationally connected subvariety in the fiber of $X'$ over $k(Z)$. Take its closure $Z'$ in $X'$. By \cite{ghs} $Z'$ is rationally connected and is clearly contained in $X'_0$. 
\end{proof}

\begin{proof}[Proof of (\ref{deg-ran})]
 Let $\{\Delta_i\}$ be the set of irreducible components of $X_0$. Let $\Delta = \sum_i \Delta_i$. We now want to 
run a relative minimal model program for the pair $(X,\Delta)$ over $C$ and apply the main theorem (see (\ref{mmp})) of \cite{mmpgen}. However, $(X,\Delta)$ is not $klt$ but only log canonical.
Nevertheless, for sufficiently small and positive $\epsilon$, $(X,\Delta-\epsilon\pi^*(0))$ 
is $klt$. And for any curve class $\eta$ in $N_1(X/C)$, $\pi^*(0)\cdot \eta=0$, running the relative minimal model program for $(X,\Delta-\epsilon\pi^*(0))$, it is equivalent to running the relative minimal model program for $(X,\Delta)$. Thus in order to run the minimal model program resulting in a Mori fibration we only need to check that $K_X+\Delta$ is not $\pi$-pseudo-effective. But in our case the generic fibre, $X_g$, of $X\to C$ is a smooth rationally connected variety. The intersection number of $K_{X_g}$ with any free rational curve in $X_g$ is negative. Thus $K_{X_g}$ is not pseudo-effective. Since $\Delta$ is vertical, this implies that $K_X+\Delta$ restricted to the generic fiber is not pseudo-effective. In particular $K_X+\Delta$ itself cannot be pseudo-effective. After running the program one arrives at the following situation.

$$\xymatrix{
(X,\Delta) \ar@{-->}[r]^{\phi}\ar[d]_{\pi} & (X',\Delta')\ar[d]^{\pi'}\ar[r]^h & Y \ar[dl]^f \\
C \ar@{=}[r]                       & C                                 &
}$$ where
\begin{enumerate}
\item[(i)] Relative Picard number of $h$ is one,
\item[(ii)] $X$ and $X'$ are birational, and
\item[(iii)] $h: X'\to Y$ is a Fano contraction, i.e. the generic fiber of $h$ is a log Fano variety.
\end{enumerate}

\noindent The generic fiber of $Y\to C$ is also a rationally connected since it is dominated by the generic fiber of $X' \to C$. Thus by induction on dimension, $Y_0=f^{-1}(0)$ has a rationally connected subvariety $Z/k$. Let $p$ be the generic point of $Z$. By (\ref{fibrationsing}), $Y$ is $klt$. The generic fiber of $X'\to Y$ is a log Fano variety. Thus by (\ref{mainthm}) there exists a rationally connected subvariety of $h^{-1}(p)$ defined over $k(p)$. Let $Z'$ be its closure in $X'$. By \cite{ghs}, $Z'$ is a rationally connected subvariety of $X'_0$. Since both $X'$ and $X$ are $klt$, an application of lemma(\ref{birational}) completes the proof. 

\end{proof}

\section{Fano Fibrations}\label{section:surface}
\noindent In this section we prove some results (of independent interests) on Fano fibrations obtained by contraction of an extremal ray. These results will then be used to prove (\ref{surface}).  

\noindent Let $\pi:(X,\Delta)\to S$ be a relative log canonical pair and with 
$X$ $\Q$-factorial with $\Delta$ effective. 
Let $h:X\to Y$ be a Fano fibration obtained by 
contracting an extremal ray. $Y$ is automatically 
$\Q$-factorial by (\cite{kollarmori}, $3.36$). Let $\Delta=\sum_ia_i\Delta_i$  where $\Delta_i$ are irreducible. 
Let $n=dim(X)$ and $r=dim(Y)$.

\begin{lem}\label{equicodim}
Let $D$ be any prime Weil divisor on $X$. Then either $h(D)= Y$ or $h(D)$ is a divisor.
\end{lem}
\begin{proof}
Assume that $\dim(h(D))=k\leq r-2$.
After cutting $Y$ and then cutting $X$ by general very ample hypersurfaces,
we may assume $h^*(H_Y)^k\cdot H_{X}^{n-k-2}$ is a surface whose
image in $Y$ is still 2 dimensional (because $r-k\geq 2$), and
$D\cdot h^*(H_Y)^{k}\cdot H_{X}^{n-k-2}$ is a curve which is
contracted to a point on $Y$, then we have $D^2\cdot
h^*(H_Y)^k\cdot H_{X}^{n-k-2} < 0$. On the other hand, let $E$
be a curve on the fibre over a general point. Since $E$ and $D$ do
not intersect $E\cdot D=0$, so $D\cdot h^*(H_Y)^{k}\cdot
H_{X}^{n-k-2}$ and $E$ cannot be numerically proportional. But this
contradicts the fact that the relative Picard number of $h$ is
one.
\end{proof}
\begin{lem}\label{oneone} Let $D$ be a prime Weil divisor on $X$ which does not
dominate $Y$. 
Then $\supp(D)=\supp(h^{-1}h(D))$.
In particular if $\Delta_i$ and $\Delta_j$ do not dominate $Y$,
then $h(\Delta_i)$ and $h(\Delta_j)$ are distinct divisors on $Y$ unless $i=j$.
\end{lem}
\begin{proof} Assume there exists an irreducible component $D'$ of $h^{-1}h(D)$ 
different from $D$. Then $h(D')\subset h(D)$. But by the previous lemma(\ref{equicodim}), 
both are irreducible divisors in $Y$. Hence $h(D')=h(D)$.
Cut $Y$ and $X$ by general very ample 
hypersurfaces so that $h^*(H_Y)^{r-1}\cdot H_{X}^{n-r-1}$ is a
surface whose image is a curve. $D'\cdot h^*(H_Y)^{r-1}\cdot H_{X}^{n-r-1}$ and 
$D\cdot h^*(H_Y)^{r-1}\cdot H_{X}^{n-r-1}$ are curves on
$h^*(H_Y)^{r-1}\cdot H_{X}^{n-r-1}$, which are mapped to the same
point on $Y$. Then by Zariski's lemma (\cite{bpv}, {\rm III}$.8.2$) they cannot be numerical
proportional which is a contradiction since $h$ has relative Picard number one. 
\end{proof}

\noindent The hypothesis of (\ref{deg-ran}) says that $X$ is smooth and the divisor $\Delta$
whose support is $X_0$ is $snc$ (each irreducible component of $\Delta$ appearing
with coefficient one).
However, we know $(X,\Delta)$ is $dlt$ (cf. \cite{kollarmori}), which is a property preserved under Log Minimal Model Program. Now we pinpoint the property we need for the $dlt$ condition as follows,
\begin{point}\label{hypommp}(Hypothesis) Let $\pi:(X,\Delta)\to S$ be 
projective log canonical pair where $\Delta$ 
is a vertical divisor. Assume $X$ is $\Q$-factorial. 
Let $\{\Delta_i\}_{1 \leq i \leq m}$ be the set of all irreducible components of 
$\Delta$ and assume each $\Delta$ appears with coefficient one in $\Delta$. 
Assume that the only log canonical centers of $(X,\Delta)$ are given by $\Delta_I$ 
where $I\subset \{1,..,m\}$ and $\Delta_I=\Intersection_{i\in I}\Delta_i$. 
\end{point}

\begin{rem}
(\ref{hypommp}) in particular implies that for each $i\in \{1,..,m\}$, $(X,\Delta_i)$ is $plt$.
Hence each irreducible component $\Delta_i$ is $klt$ by inversion of adjunction. 
We will show (see (\ref{dltklt})) that even after performing divisorial contractions and 
flips, the resulting divisor has $klt$ components. The following proposition proves
the same for Fano fibrations. 
\end{rem}

\begin{prop}\label{fibrationsing}
Let $(X,\Delta)$ satisfy hypothesis(\ref{hypommp}). Let $X\to Y$ be a Fano fibration
obtained by contracting an extremal ray. Let $D_i$ be the image of $\Delta_i$ in $Y$ and
let $D= \sum_{i=1}^mD_i$. 
Then $(Y,D)$ also satisfies (\ref{hypommp}).
\end{prop}

\begin{rem}To prove (\ref{fibrationsing}) we need Kawamata's canonical bundle formula 
as proved in \cite{Kol05}. A similar method  
has been used before in \cite{ambro} to study singularities of Fano contractions.
\end{rem}

\begin{defn}\cite{Kol05}(Standard normal crossing assumptions)\label{nc}
We say that a projective morphism 
$f:X\to Y$ together with $\Q$-divisors $R$, $B$ satisfy the standard normal crossing assumptions 
if the following hold
\begin{enumerate}
\item $X$,$Y$ are smooth,
\item $R+f^*B$ and $B$ are $snc$ divisors,
\item $f$ is smooth over $Y\setminus B$, and
\item $R$ is a relative $snc$ divisor over $Y\setminus B$.
\end{enumerate}
\end{defn}

\begin{thm}[Kawamata's canonical bundle formula, see \cite{Kol05}]\label{KK}
Let $f:X\to Y$ and $R$,$B$ satisfy the standard normal crossing assumptions. 
Let $F$ be the generic fibre of $f$. Assume that $K_X+R\sim_{\Q} f^*H_R$ for some
$\Q$-divisor $H_R$ on $Y$. Let $R=R_h+R_v$ be the horizontal and vertical parts of $R$ and assume that when we write
$R_h=R_{h(\geq 0)}-R_{h(\leq 0)}$ as the difference of its positive and negative parts, we have
$h^0(F,\rdown{R_{h(\leq 0)}|_F})=1$. Then we can write
$$K_{X/Y}+R \sim_\Q f^*(J(X/Y,R)+B_R)$$
where
\begin{enumerate}
\item The moduli part $J(X/Y,R)$ is nef and depends only on $(F,R_h|_F)$ and $Y$,
\item The boundary part $B_R$ depends only on $f:X \to Y$ and $R_v$. More precisely, $B_R$ is the unique smallest $\Q$-divisor supported on $B$ such that
$$red(f^*B)\geq R_v+f^*(B-B_R).$$
Moreover,
\item The pair $(Y,B_R)$ is $lc$ iff The pair $(X,R)$ is $lc$.
\item If $\rdown{R_h}\leq 0$, then the pair $(Y,B_R)$ is $klt$ iff the pair $(X,R)$ is $klt$.
\item $B_i$ appears in $B_R$ with nonnegative coefficient if $f_*\sO_X(\rdown{-R_v})\subset \sO_{B_i,Y}$
\end{enumerate}
\end{thm}
\begin{rem}
In the above formula, $J(X/Y,R)$ is only defined as a $\Q$-divisor class. However, $B_R$ is a $\Q$-divisor.
\end{rem}

\begin{proof}[Proof of (\ref{fibrationsing})]
Since $h$ is a contraction of a $-(K_X+\Delta)$-negative ray, $-(K_X+\Delta)$ 
is $h$-ample. Thus there exist ample $\Q$-divisors $H_X$ (resp. $H_Y$) on $X$ (resp. $Y$)
in general position such that 
$$ K_X + \Delta + H_X \sim_{\Q} h^*(H_Y)$$
Choose log resolutions $g:\overline{X} \to (X,\Delta)$ and $d:\overline{Y}\to Y$ such 
that we have the following commutative diagram
$$\xymatrix{
\overline{X} \ar@{->}[r]^{g}\ar[d]_{\overline{h}} & X\ar[d]^{h} \\
\overline{Y}\ar@{->}[r]^{d}                                   & Y
}$$
There exists a Cartier divisor $R$ on $\overline{X}$ such that 
$$ K_{\overline{X}/X}+R= g^*(\Delta+H_X)$$
After further blowup we may assume $\overline{h}:\overline{X}\to \overline{Y}$, the 
divisor $R$, and a reduced divisor $B$ on $\overline{Y}$ containing $\supp(d^{-1}(D))$
satisfy the standard normal crossing assumptions (\ref{nc}). 

\noindent We can choose above $H_X=\sum d_iH_i$, where $0<d_i \ll 1$ and $H_i$ are divisors which are in general position.   We now apply (\ref{KK}) to $\overline{h}:\overline{X}\to \overline{Y},R$ and $B$. Since the negative part of $R$ consists of 
exceptional divisors,  $h^0(F,\rdown{R_{h(\leq 0)}})=1$. Thus by (\ref{KK})
we have 
$$K_{\overline{X}}+R \sim_{\Q} 
\overline{h}^*(K_{\overline{Y}}+J(\overline{X}/\overline{Y},R)+B_R)$$
\begin{lemma}
For every $i$ the coefficient in $B_R$ of the birational
transform of $D_i$ in $\overline{Y}$ is $1$. Furthermore, every log canonical center $W$ of
$(\overline{Y},B_R)$ satisfies that $d(W)=\intersection_J D_j$ for
some $J\subset \{1,...,m\}$.
\end{lemma}
\begin{proof}
The first part of the lemma is because of (\ref{KK}(2)) and the fact that $R$ is a subboundary (i.e. the coefficients are less or equal to 1, not necessarily nonnegative) which contains a 
component dominating $D_i$ for every $i$, whose coefficient in $R$ is 1. 

\noindent The reason for the second part is similar to the
one just given. Let $D$ be an exceptional divisor of $\overline{Y}$ 
whose center in $\overline{Y}$ is a log canonical center $W$
of $(\overline{Y},D)$. To show $d(W)$ is of the form $\Intersection_JD_J$
in $Y$, we may blow up $\overline{Y}$ and $\overline{X}$ further 
and assume $D=W$ and
the coefficient of $D$ in $B_R$ is $1$. But this implies that there is a 
component $E$ of $R_v$ dominating $B$ such that the coefficient of 
$E$ in $R_v$ is one. The image of $E$ in $X$ 
is a log canonical center of $(X,\Delta)$ and hence is of the form
$\Intersection_{i\in J}\Delta_i$ for some $J\in \{1,...,m\}$.
Thus $d(W)=\Intersection_{i\in J}D_J$. 
\end{proof}

\noindent Now we only need to prove that $(Y,D)$ is log canonical. For this we first
claim that
\begin{lemma}
$d_*B_R$ is an effective divisor.
\end{lemma}
\begin{proof}
By (\ref{KK}($2$)) it follows 
that if an irreducible component $B_1$ appears in $B_R$ with negative coefficient
then every component of $\overline{h}^*B_1$ appears in $R_v$ with negative coefficient.
But since all negative components of $R$ have to be $g$ exceptional, it follows
that $B_1$ is also $d$-exceptional. This proves the lemma. 
\end{proof}
\noindent To check $(Y,D)$ is canonical, it is enough to check that for every divisor
$E$ on $\overline{Y}$ $a(E,Y,D)\geq -1$. Since $J(\overline{X}/\overline{Y})$ 
is nef, for every $\epsilon >0$ we can choose an effective divisor 
$J_{\epsilon}\sim_{\Q} J(\overline{X}/\overline{Y})+d^*H_Y$ such that 
for every divisor $E$ on $Y$ $a(E,\overline{Y},J_{\epsilon}+B_R)>-1-\epsilon$.
Because $K_{\overline{Y}}+J_{\epsilon}+B_R \sim_{\Q} d^*(A)$ for some
$\Q$-divisor $A$ on $Y$, then we have
$$K_{\overline{Y}}+J_{\epsilon}+B_R \sim_{\Q} d^*(K_Y+d_*(J_{\epsilon}+B_R))$$
But $D\leq d_*(J_{\epsilon}+B_R)$. Hence for every divisor $E$ on $Y$ we have 
$$ a(E,Y,D)\geq -1-\epsilon $$
Taking the limit as $\epsilon \to 0$ we see that $(Y,D)$ is log canonical.
\end{proof}

\begin{proof}[Proof of (\ref{surface})]
The idea of the proof is similar to that of (\ref{deg-ran}). Let us assume $\dim X_t=3$, for the case of smaller dimension, the argument is similar but eaiser. Let $\Delta$ be the sum of all irreducible components of $X_0$, each with coefficient one. Then after running the relative minimal model program for $(X,\Delta)$ we arrive at one of the following two cases.
\begin{enumerate}
\item[Case(1)]
$$\xymatrix{
(X,\Delta) \ar@{-->}[r]^{\phi}\ar[d]_{\pi} & (X',\Delta')\ar[d]^{\pi'} \\
C \ar@{=}[r]                       & C
}$$ where $\pi'$ has relative Picard number one.
\item[Case(2)]
$$\xymatrix{
(X,\Delta) \ar@{-->}[r]^{\phi}\ar[d]_{\pi} & (X',\Delta')\ar[d]^{\pi'}\ar[r]^h & Y \ar[dl]^f \\
C \ar@{=}[r]                       & C                                 &
}$$ where
\begin{enumerate}
\item[(i)] Relative Picard number of $h$ is one,
\item[(ii)] $X$ and $X'$ are birational, and $h: X'\to Y$ is a Fano contraction of relative dimension two, and
\item[(iii)] $Y/C$ is of relative dimension one with generic fiber isomorphic to $\P^1_{k(C)}$.
\end{enumerate}
\item[Case(3)]
$$\xymatrix{
(X,\Delta) \ar@{-->}[r]^{\phi}\ar[d]_{\pi} & (X',\Delta')\ar[d]^{\pi'}\ar[r]^h & Y \ar[dl]^f \\
C \ar@{=}[r]                       & C                                 &
}$$ where
\begin{enumerate}
\item[(i)] Relative Picard number of $h$ is one,
\item[(ii)] $X$ and $X'$ are birational, and $h: X'\to Y$ is a Fano contraction of relative dimension one, and
\item[(iii)] $Y/C$ is of relative dimension two.
\end{enumerate}
\end{enumerate}

\noindent {\it Case ($1$)}: Every irreducible component of $\Delta'$ contributes to the relative Picard number of $\pi'$. Since the relative Picard number is one $\Delta'$ must be irreducible. It is also geometrically connected since fibres of $\pi'$ are geometrically connected. By lemma(\ref{dltklt}) $\Delta'$ is $klt$, in particular normal. A normal geometrically connected variety over $k$ is geometrically irreducible. Hence $\Delta'$ is geometrically irreducible. 
\noindent  Since $\Delta'$ is a degeneration of a proper rationally chain connected scheme, it is rationally chain connected.  But $\Delta'$ is $klt$, hence by \cite{hm05} it is in fact rationally connected. This proves that $\Delta'$ is geometrically irreducible as well as rationally connected. By (\ref{dltklt}($i$)), $\Delta' = \phi_*\Delta_i$ for some irreducible component $\Delta_i$ of $\Delta$. Hence $\Delta_i$ is also rationally connected and geometrically irreducible.\\

\noindent {\it Case ($2$)}:  
Let $\{\Delta_i'\}_{1\leq i \leq s}$ be the irreducible components of $\Delta'$ and let $\Delta_i'=\phi_*\Delta_i$.
Let $D=\sum_{i=1}^sD_i$. By (\ref{fibrationsing}) $(Y,D)\to C$ satisfies hypothesis (\ref{hypommp}). 

\noindent We first claim that there exists an irreducible components $D_{j_0}$ of $D$ which is rationally connected and in particular geometrically irreducible. The generic fibre of $Y/C$ is isomorphic to $\P^1$ and hence $-(K_Y+D)$ is not pseudo-effective. Thus after running the minimal model program for $(Y,D)/C$ we arrive at the following situation. 
$$\xymatrix{
(Y,D) \ar@{-->}[r]^{\phi}\ar[d]_{\pi_Y} & (Y',D')\ar[d]^{\pi'_Y} \\
C \ar@{=}[r]                       & C
}$$ where $\pi'_Y$ has relative Picard number one. The claim now follows by an argument similar to that in Case($1$). 

\noindent Now look at $\Delta'_{j_0}$. It is geometrically irreducible by (\ref{oneone}). Let $\Delta''\to \Delta_{j_0}$ be a resolution of singularities. Then by \cite{hm05}(1.2), every fibre of $h':\Delta''\to D_{j_0}$ is rationally chain connected. In particular the generic fiber is rationally connected since it is also smooth. This together with \cite{ghs} and the fact that $D_{j_0}$ is rationally connected proves that $\Delta''$ and hence $\Delta'_{j_0}$ is rationally connected. Since $\Delta_{j_0}$ and $\Delta_{j_0}'$ are birational, $\Delta_{j_0}$ is also rationally connected.\\

\noindent{\it Case ($3$)} The proof is similar as in Case (2). We need to use induction to prove there exists a rationally connected subvariety in $Y_0$. Hence it suffices to prove that when we run the minimal model program for $(Y,D)$, it terminates with a Fano contraction. In another words, it suffices to show the non-pseudo-effectivity of $K_Y+D$. Now $\phi|_{X_t}: X_t\to Y_t$ is a Fano contraction of an extremal ray of $K_{X_t}$, and $X_t$ is a terminal three-fold. Apply the Theorem 1.2.7 of \cite{mp}, we conclude that $Y_t$ only contains Du Val singularity of type A. In particular, since $X$ is rationally connected, $Y_t$ is a rational surface with canonical singularity. Take the minimal resolution $d:\overline{Y}_t \to Y_t$, then
$$H^0(Y_t,mK_{Y_t})=H^0(\overline{Y}_t,d^*(mK_{Y_t}))=H^0(\overline{Y}_t,mK_{\overline{Y}_t})=0.$$
Since $K_Y+D|_{Y_t}=K_{Y_t}$, we conclude that $K_Y+D$ is non-pseudo-effective.   
\end{proof}

\begin{rem}\label{hard}

\noindent One can also try to answer (\ref{mainq'}) or (\ref{mainq}) in higher dimensions by the same method as above. In this case one would not have to deal with the case of degeneration of Fano varieties separately. The missing property here is the non-pseudo-effectivity of $K_Y+D$ where $(Y,D)$ is a result of a Fano contraction. This non-pseudo-effectivity is by no means automatic except in the case when relative dimension of $Y/C$ is one. In fact, in dimension $\geq 2$, Koll\'ar has constructed a family of rational varieties with only quotient singularities, but whose canonical divisor is $\Q$-ample (cf. \cite{kollar06}). 
As a stronger evidence to show the restriction of our approach toward (\ref{mainq}), Prokohorv even constructs an example of $h:X\to Y$ a Fano contraction, where $X$ is a rationally connected variety of terminal singularity, but $Y$ has a $\Q$-effective canonical bundle.

\end{rem}
\begin{ex}[Prokhorov]\label{counter}
\noindent Let $S$ be the surface $x^n+y^n+z^n+w^n=0$ in $\P^3$, the group $G=\Z/n$ acts on $\P^3$ and $S$ with weights $(0,0,1,1)$.
Let $Y:=S/G$, and $\pi:S \to Y$ is the quotient morphism. Since $\pi$ is \'etale outside finite points, thus $\pi^*(K_Y)=K_S$. So when $n\ge 4$, $K_Y$ is $\Q$-effective. We claim $Y$ is rational. In fact, the projection map $p:S\to \P^1$ by sending $(x,y,z,w)$ to $(x,y)$ will factor through $Y$. The generic fiber of the map $p$ is birational to the affine curve $(z^n+w^n+\lambda=0)\subset \A^2$, so the generic fiber of the map from $Y$ to $\P^1$ is birational to the curve $(z^n+w^n+\lambda=0)/ G$, which is also rational. We conclude that $Y$ is a rational surface.

\noindent Now consider the action of $G(=\Z/n)$ on $\P^{n-1}$ with the weights $(0,1,...,n-1)$. We define $X:= (\P^{n-1}\times S) /G$, where $G$ acts diagonally. For a fixed point $x$ of a $g\in G$ on $(\P^{n-1}\times S) /G$, if we write the eigenvalue of $g$ on $T_x$ as $e(r_1),...,e(r_{n+1})$, where $e(x):=e^{2\pi ix}$ and $0\le r_i<1$. Then to prove $X$ contains only terminal singularity, it suffices to prove that for any element nontrivial element $g$ and any $x$ stablized by $g$, the $age$ satisfies
$$\age_x(g):=r_1+\cdots+r_{n+1}>1.$$ 
(cf. \cite{reid}, Theorem 4.6). Now when $n\ge 4$, we can see that the above condition holds for our case.

\noindent We define the map $h:X\to Y$ by descending the second projection of $p_2: \P^{n+1}\times S\to S$, which is equivariant. It is easy to see that $h$ is a Mori contraction of an extremal ray. Since any terminal variety is an output of a minimal model program which starts with a smooth variety (simply reverse the desingularization process), we in fact have an example, where we start with a smooth rationally connected variety, but running the minimal model program terminates with a contraction to a variety of ample canonical bundle. 
\end{ex}

\noindent Finally we give an example of a simple normal crossing variety of dimension $3$ which is rationally chain connected but none of its components are rationally connected. 
\begin{ex}\label{counterex}
\noindent Let $C$ be any smooth curve in $\P^2$ of genus $g \geq 2$. Let
$i:C \inj \P^2$ be the closed embedding of $C$ in $\P^2$.
Let
$$
\xymatrix{
            &  \P^2\times_k C = X_1 \\
C\times_k C\ar[ru]^{i_1}\ar[rd]^{i_2} &   \\
             &  \P^2\times_k C = X_2
}
$$
be two maps defined by
$i_1 = (i,id_{C})$ and $i_2 = (id_{C}, i)$. Glue $X_1$ and $X_2$ along the maps $i_1,i_2$. The
resulting three dimensional variety $X=X_1\Union X_2$ is simple normal crossing as well as
rationally chain connected. But neither $X_1$ nor $X_2$ is rationally connected.
\end{ex}


\section{Rationally Connected Subvarieties of the Exceptional Locus}\label{section:exc}

\noindent In this section we prove Theorem(\ref{exc.thm}).

\begin{lem}
Let $X/k$ be a normal variety. Let $\pi:X'\to X$ be a proper birational morphism.  
Then after further blowing up $X'$ if necessary we can find an effective exceptional
divisor $A$ on $X'$ such that $-A$ is $\pi$-ample.
\end{lem}
\begin{proof}
It is enough to prove the lemma for one particular $\pi:X'\to X$ with $X'$ smooth. 
For any quasiprojective variety $S$, if $f:S_{\sI} \to S$ is the blow up of an ideal sheaf $\sI$ 
on $S$ then $E(\sI)= \pi^{-1}\sI$ is an effective divisor on $S_{\sI}$ which is anti $f$-ample.
Now let $\pi:X'\to X$ be the resolution of singularities obtained by successively 
blowing up centers disjoint from the smooth locus of $X$, in particular $\pi$ is projective. Suitable positive 
linear combination of (strict transforms of) divisors of the form $E(\sI)$ gives
us an effective divisor $A$ on $X'$ such that $-A$ is $\pi$-ample. Since $X$ 
is normal this divisor $E$ is $\pi$-exceptional.
\end{proof}

\begin{proof}[Proof of (\ref{exc.thm})]
For getting hold of a geometrically irreducible subvariety of $\pi^{-1}(p)$ 
we follow the general idea of the proof of main theorem in \cite{kollarax}.
First by base extending to $k(p)$ we may assume $p$ is a $k$-point.
By replacing $X$ by a suitable open neighbourhood of $p$ we may assume $X$ is 
affine.
We claim that there exists an effective $\Q$-divisor 
$\Delta_1\sim_\Q 0$ passing through $p$ such that 
$(X,\Delta+\Delta_1)$ is $klt$ outside $p$ and log canonical at $p$.
Since $X$ is affine, taking a $\Delta_1$ to be a suitable rational multiple of 
a sufficiently general section of $\sO_X$ passing through $p$ with high multiplicity
does the job. \\

\noindent By further blowing up $X'$ we may assume 
\begin{enumerate}
\item[(i)]$\pi:X'\to X$ is a log resolution of $(X,\Delta+\Delta_1)$,
\item[(ii)]$\pi^{-1}(p)$ is a divisor, and
\item[(iii)]There is an effective exceptional divisor $A$ on $X'$ such that 
$-A$ is $\pi$-ample.
\end{enumerate}
Let $\{E_i\}_{1\leq i \leq r}$ be the exceptional divisors lying over $p$ and let 
$\{F_j\}_{1\leq j \leq s}$ be the remaining exceptional divisors of $\pi$.\\

\noindent Henceforth we will use the notation "$\sim_\Q $" to denote 
$\Q$-numerical equivalence of $\Q$-divisors relative to $\pi$. In another words, for two $\Q$-Cartier divisors $A$ and $B$, we write $A\sim_{\Q}B$ iff
$$A\cdot \eta=B\cdot \eta,$$
for any $\eta \in N_1(X'/X)$, then we have
$$ K_{X'} \sim_{\Q} \sum_{i=1}^r a_iE_i + \sum_{j=1}^sb_jF_j - \Delta'$$ 
where 
\begin{enumerate}
\item[(i)] $a_i\geq -1  \ \forall \ i$ with equality holding for at least one $i$,
\item[(ii)] $b_j >-1\ \forall j$, and
\item[(iii)] $\Delta'$ is the birational transform of $\Delta+\Delta_1$ hence $\rdown{\Delta'}=0$.
\end{enumerate}
Let 
$$ A = \sum_{i=1}^rc_i E_i + \sum_{j=1}^sd_jF_j$$
By replacing $A$ by $\frac{1}{n}A$ for $n$ large enough we may assume the 
coefficients $c_i$ and $d_j$ are sufficiently small as compared to one. 
For small rational numbers $\epsilon>0$ and $\{\delta_i\}_{1\leq i \leq r}$
consider the following divisor 
$$K_{X'}-\pi^*(K_X+\Delta+(1-\epsilon)\Delta_1)-(A+\sum_i\delta_iE_i)  \sim_\Q  K_{X'}-(A+\sum_i\delta_iE_i) -\Delta'(\epsilon)$$
$$\sim_{\Q}\sum_{i=1}^r a_i(\epsilon,\delta)E_i + \sum_{j=1}^sb_j(\epsilon)F_j - \Delta'(\epsilon) $$
where $\Delta'(\epsilon)$ is the birational transform of $\Delta+(1-\epsilon)\Delta_1$. 
Notice that $b_j(\epsilon)>-1$ for all $\epsilon >0$ and all $j$. 
Choose $\epsilon$ and $\delta_i's$ such that 
\begin{enumerate}
\item[(i)] $A+\sum_i\delta_iE_i=A'$ is effective and such that $-A'$ is $\pi$-ample
(this is true for all $\delta_i$ small enough because ampleness is an open condition), and
\item[(ii)]$a_i(\epsilon,\delta)\geq -1$ with equality holding for exactly one $i$ (say $i=1$).
\end{enumerate}
Thus we can write 
$$ K_{X'}-A' \sim_{\Q} -E_1 + E - M$$
where  
$$ E = \sum_{i\geq 2}\rup{a_i(\epsilon,\delta)}E_i+\sum_j\rup{b_j(\epsilon)}F_j$$
is an effective exceptional divisor and $M$ is an effective $snc$ divisor (having no 
common components with $E_1$) and such that $\rdown{M}=0$.

\noindent Thus by relative Kawamata-Viehweg vanishing theorem we get 
$$ R^1\pi_*\sO_{X'}(-E_1+E)=0.$$
This gives us a surjection 
$$\sO_X \cong \pi_*\sO_{X'}(E)\to \pi_*\sO_{E_1}\to 0$$
which proves $E_1$ is geometrically irreducible. \\

\noindent We claim that $E_1$ is also rationally connected. 
Since $-A'$ is $\pi$-ample and since $X$ is affine it is actually ample. 
Choose a general divisor $W$ such that $-A'\sim_{\Q}W$ and such that $\rdown{W}=0$.
Since $W$ is general, it has simple normal crossings with all exceptional 
divisors and also with $M$. Let $M'=M+W$. Then $M'$ supports an ample divisor and
we have 
$$K_{X'}+E_1+M' \sim_{\Q} E$$
But $E$ is an effective divisor having no common components with $E_1$.
Thus $K_{E_1}+M'_{|E_1}$ is effective and as in proof of (\ref{mainthm}) rational
connectivity of $E_1$ will follow by (\cite{hm05}, $4.1$) 
once we show $\kappa(K_{E_1}+M'_{|E_1})=0$. As in the proof of (\ref{specialcase}) 
this follows from (\ref{extensionthm}) and using the fact that $E$ is effective 
exceptional.
\end{proof}

\begin{rem}
Theorem(\ref{exc.thm}) can fail for varieties having log canonical singularities.
See (\cite{kollarax}, $9$) for an example. 
\end{rem}


\vspace{5mm}
\noindent Amit Hogadi\\
\noindent Department of Mathematics, Princeton University.\\
\noindent {\it Email:} amit@math.princeton.edu\\

\noindent Chenyang Xu\\
\noindent Department of Mathematics, Princeton University.\\
\noindent {\it Email:} chenyang@math.princeton.edu


\begin{thebibliography}{BDPP04}

\bibitem[Am04]{ambro}
Ambro, F.; 
Shokurov's boundary property.
{\it J. Differential Geom.} {\bf 67} (2004), no.{\bf 2}, 229--255. 

\bibitem[BPV84]{bpv}
Barth, W.; Peters, C.; Van de Ven, A.; Compact complex surfaces. 
Ergeb. Math. Grenz.{\bf 3}, {\bf 4}. 
Springer-Verlag, Berlin, 1984.

\bibitem[BCHM06]{mmpgen}
Birkar C.; Cascini P.; Hacon C.; McKernan J.;
Existence of minimal models for varieties of log general type,
math.AG/0610203. (2006). 


\bibitem[Ca92]{Ca}
Campana, F.;
Connexit\'e rationelle des vari\'et\'es de Fano. {\it Ann. Sci. \'Ecole Norm. Sup.} {\bf 25} (1992) 539-545.

\bibitem[dJS03]{ds03}
de Jong, A. J.; Starr, J.;
Every rationally connected variety over the function field of a curve has a rational point.
{\it Amer. J. Math.}  {\bf 125}  (2003),  no. {\bf 3}, 567--580.

\bibitem[Es03]{He}
Esnault, H.; Varieties over a finite field with trivial Chow group
of 0-cycles have a rational point.  {\it Invent. Math.} {\bf 151}
(2003), no. {\bf 1}, 187--191. 

\bibitem[EX07]{ex07}
Esnault, H.; Xu, C.; Congruence for rational points over finite fields and coniveau over local fields,
math.NT/0706.0972, (2007). To appear in {\it Trans. Amer. Math. Soc.}


\bibitem[FR05]{nr}
Fakhruddin, N.; Rajan, C.S.; Congruences for rational points on varieties over finite fields.
{\it Math. Ann.} {\bf 333}  (2005),  no.{\bf 4}, 797--809.

\bibitem[GHS03]{ghs}
Graber, T.; Harris, J.; Starr, J.;
Families of rationally connected varieties.
{\it J. Amer. Math. Soc.} {\bf 16} (2003), no.{\bf  1}, 57--67.

\bibitem[GHMS05]{mazur}
Graber, T.; Harris, J.; Mazur, B.; Starr, J.; Rational Connectivity and sections of
families over curves.  {\it Ann. Sci. \'Ecole Norm. Sup. (4)}  {\bf 38}  (2005),
no. {\bf 5}, 671--692.



\bibitem[HM06]{hm06}
Hacon, C.; McKernan, J.;
Boundedness of pluricanonical maps of varieties of general type,
{\it Invent. Math.} {\bf 166} (2006) no. {\bf 1}, 1--25. 

\bibitem[HM07]{hm05}
Hacon, C.; McKernan, J.; Shokurov's Rational Connectedness Conjecture,
{\it Duke Math. J.} {\bf 138} (2007) no. {\bf 1}, 119-136.


\bibitem[Kol96]{kollarrc}
Koll\'ar, J.; {\it Rational curves on algebraic varieties.},
Ergeb. Math. Grenz. {\bf 3}. Folge, {\bf 32}.
Springer-Verlag, Berlin, 1996.

\bibitem[Kol07a]{kollarax}
Koll\'ar, J.;  A conjecture of Ax and degenerations of Fano varieties,
{\it Israel J. Math.} {\bf  162}  (2007) no. {\bf 1}, 235--251. 

\bibitem[Kol07b]{Kol05}
Koll\'ar, J.;  Kodaira's canonical bundle formula and subadjunction, Ch. $8$ in {\it Flips for $3$-folds
and $4$-folds}, editor A. Corti, Oxford Lecture Series in Mathematics and its 
applications {\bf 35}, Oxford University Press, 2007.

\bibitem[Kol08]{kollar06}
Koll\'ar, J.;  Is there a topological Bogomolov--Miyaoka--Yau inequality.
 {\it Pure Appl. Math. Q.} {\bf 4}  (2008) no.{\bf 2},  203--236.

\bibitem[KM98]{kollarmori}
Koll\'ar, J.; Mori, S.; {\it Birational Geometry of Algebraic Varieties},
Cambridge Tract. in Math. {\bf 134}, Cambridge University Press, Cambridge, 1998.

\bibitem[KMM92a]{KMM}
Koll\'ar, J.; Miyaoka, Y.; Mori, S.;
Rational Connectedness and Boundness of Fano Manifolds,
{\it J. Diff. Geom.} {\bf 36} (1992), no. {\bf 3}, 765--779.

\bibitem[KMM92b]{kmm}
Koll\'ar, J.; Miyaoka,Y.; Mori,S.;
Rationally connected varieties.
{\it J. Algebraic Geom.} {\bf 1} (1992), no. {\bf 3}, 429--448.

\bibitem[MP06]{mp}
Mori, M.; Prokhorov, Y.; On $\Q$-conic bundles,
math.AG/0603736, (2006). To appear in  {\it Publ. Res. Inst. Math. Sci.} 


\bibitem[Re87]{reid}
Reid, M.;
Young person's guide to canonical singularities. {\it Algebraic geometry, Bowdoin, 1985 (Brunswick, Maine, 1985)}, 345--414,
Proc. Sympos. Pure Math., 46, Part 1, {\it Amer. Math. Soc., Providence, RI}, 1987. 

\bibitem[St06]{starr}
Starr, J.; Degenerations of rationally connected varieties and PAC fields,
math.AG/0602649, (2006).

\bibitem[Zh06]{zhang}
Zhang, Q.; Rational connectedness of log $\Q$-Fano varieties.
{\it J. Reine Angew. Math.} {\bf 590} (2006), 131--142.


\end{thebibliography}
\end{document}